\documentclass[a4paper,reqno]{amsart}

\usepackage{amsmath}
\usepackage{amsthm}
\usepackage{amsfonts}
\usepackage{amssymb}

\usepackage{hyperref}

\hypersetup{
    colorlinks=false,
    pdfborder={0 0 0},
    pdfstartview={XYZ null null 1.00}
}

\numberwithin{equation}{section}

  \def\<{\langle}
  \def\>{\rangle}

  \def\ker{\mathrm{Ker}\,}
  
  \def\inv{\mathrm{Inv}\,}
  \def\inter{\mathrm{int}\,}

  \def\o{\overline}
  \def\h{\widehat}

  \newcommand{\x}{\bar u}
  \newcommand{\y}{\bar v}
  \newcommand{\z}{\bar w}

  \def\d{ \, d }

  \def\R{\mathbb{R}}

  \def\C{\mathbb{C}}

  \def\ffor{ \qquad \mathrm{for} \quad }
  \def\aas{ \qquad \mathrm{as} \quad }

\theoremstyle{plain}
  \newtheorem{theorem}{Theorem}[section]

  \newtheorem{lemma}[theorem]{Lemma}

\theoremstyle{definition}

  \newtheorem{remark}[theorem]{Remark}

\begin{document}

\title[Connecting orbits for nonlinear...]{Connecting orbits for nonlinear differential equations at resonance}

\author{Piotr Kokocki}
\address{\noindent Faculty of Mathematics and Computer Science \newline Nicolaus Copernicus University \newline Chopina 12/18, 87-100 Toru\'n, Poland}
\email{pkokocki@mat.umk.pl}
\thanks{The researches supported by the MNISzW Grant no. N N201 395137}

% \subjclass[2010]{47J35, 35B10, 37L05}


\begin{abstract}
We study the existence of orbits connecting stationary points for the first order differential equations being at resonance at infinity, where the right hand side is the perturbations of a sectorial operator. Our aim is to prove an index formula expressing the Conley index of associated semiflow with respect to appropriately large ball, in terms of special geometrical assumptions imposed on the nonlinearity. We also prove that the geometrical assumptions are generalization of well known in literature Landesman-Lazer and strong resonance conditions. Obtained index formula will be used to derive the criteria determining the existence of orbits connecting stationary points for the heat equation being at resonance at infinity.
\end{abstract}

\maketitle

\setcounter{tocdepth}{2}
%\tableofcontents

\section{Introduction}

We consider nonlinear differential equations of the form
\begin{align}\label{row-par}
\dot u(t)  = - A u(t) + \lambda u(t) + F (u(t)),  \qquad  t > 0
\end{align}
where $\lambda$ is a real number, $A: X\supset D(A)\to X$ is a sectorial operator on a Banach space $X$ and $F:X^\alpha\to X$ is a continuous map. Here $X^\alpha$ for $\alpha\in(0,1)$, is a fractional power space associated with $A$. We intend to study the existence of orbits connecting stationary points for the equation \eqref{row-par} being at {\em resonance at infinity}, that is, $\ker(\lambda I - A)\neq \{0\}$ and $F$ is a bounded map. To explain this more precisely, assume that, for every initial data $x\in X^\alpha$, the equation \eqref{row-par} admits a (mild) solution $u(\,\cdot\,; x):[0,+\infty)\to X^\alpha$ starting at $x$. We can define \emph{the semiflow} $\Phi:[0,+\infty)\times X^\alpha \to X^\alpha$ given by
\begin{equation*}
\begin{aligned}
\Phi(t, x):= u(t;x) \ffor t\in[0,+\infty), \ x\in X^\alpha.
\end{aligned}
\end{equation*}
Then {\em the stationary point} of \eqref{row-par} is an element $x_0\in X^\alpha$ such that $\Phi(t,x_0) = x_0$ for $t\ge 0$ and {\em the full solution} (or {\em orbit}) of this equation is a map $u:\R\to X^\alpha$ that satisfies the following equality
\begin{equation*}
    \Phi(s,u(t)) = u(t + s) \ffor s\ge 0, \ t\in\R.
\end{equation*}
We say that the full solution $u$ connects stationary points $u_+,u_-\in X^\alpha$ provided there are sequences $(t^+_n)$ and $(t^-_n)$ with $t^\pm_n\to \pm\infty$ as $n\to +\infty$ such that
$$\lim_{n\to +\infty} u(t^+_n) = u_+ \qquad\text{and}\qquad \lim_{n\to +\infty} u(t^-_n) = u_-.$$
A tool that we will use to study this problem is a version of the Conley index for semiflows defined on metric spaces. This index was introduced initially for semiflows acting on finite dimensional vector spaces, see \cite{MR511133}, \cite{MR797044}. In \cite{MR637695} and \cite{MR910097}, Rybakowski extended this index theory on arbitrary metric space, which gave a rise to study the dynamics of partial differential equations.

This paper is motivated by \cite{MR731150}, \cite{MR798176} and \cite{MR1992823} where the Conley index methods were used to prove that existence of orbits connecting stationary points for the equation \eqref{row-par}. However in these articles authors studied this problem assuming lack of resonance at infinity. In this paper, we continue to study this problem in the case of resonance at infinity. The main difficulty lies in the fact that, in the presence of resonance at infinity, the problem of existence of orbits connecting stationary points may not have solution for general nonlinearity $F$. This fact has been explained in detail in Remark \ref{rem-res}. Therefore, our aim is to overcome this problem by new theorems determining the existence of connecting orbits for \eqref{row-par}, in terms of appropriate geometrical assumptions imposed on $F$.

\noindent In Section 2 we briefly recall a necessary properties of Conley index that we will use in the paper, such as existence, multiplicity and homotopy invariance.

\noindent In Section 3 we formulate geometrical assumptions $(G1)$ and $(G2)$ (see page \pageref{g1g2}) and use them to prove our main result, {\em the index formula for bounded orbits}, which express the Conley index of the invariant set of $\Phi$ in sufficiently large ball, in terms of assumptions $(G1)$ and $(G2)$.

\noindent Finally, in Section 4 we provide applications for particular partial differential equations. First of all, in Theorems \ref{lem-est2} and \ref{lem-est3}, we prove that if $F$ is a Niemytzki operator associated with a map $f:\Omega\times\R\to\R$, then the well known in literature Landesman-Lazer and strong resonance conditions are actually particular case of assumptions $(G1)$ and $(G2)$. The Landesman-Lazer conditions were introduced in \cite{MR0267269} to study the stationary points of heat equations being at resonance at infinity. There is many papers using these conditions to study the heat or telegrapher equations in the case of resonance. For instance in \cite{MR0513090} and \cite{MR0487001} theorems for the existence periodic solutions or stationary points were derived. Similarly, the strong resonance conditions were studied for example in \cite{MR713209}, \cite{MR597281} in order to obtain the existence of stationary points and periodic orbits of heat equation. Encouraged by these results we use the abstract results obtained in Section 3 and prove effective criteria determining the existence of orbits connecting stationary points for the nonlinear heat equation, in terms of Landesman-Lazer and strong resonance conditions.

\section{Homotopy index}

In this section we gather the properties of Conley Index which are necessary in this paper. Fore more details see \cite{MR910097}. The continuous map $\Phi:[0,+\infty)\times X^\alpha \to X^\alpha$ is a semiflow on $X^\alpha$ provided $$\Phi(0,x) = x \ \ \text{and} \ \ \Phi(t+t',x) = \Phi(t, \Phi(t',x)) \ffor t,t'\ge 0, \ x\in X^\alpha.$$
A map $\sigma:[-\delta_1,\delta_2) \to X^\alpha$, where $\delta_2 > 0$, $\delta_1 \ge 0$, is said to be  \emph{a solution} of $\Phi$, if $$\Phi(t,\sigma(s)) = \sigma(t+s) \ffor t\ge 0 \text{ and } s\in [-\delta_1,\delta_2), \ t + s\in[-\delta_1,\delta_2).$$
If $\sigma$ is defined on $\R$, then $\sigma$ is called \emph{a full solution}. We say that the full solution $\sigma$ {\em connects} stationary points $u_+,u_-\in X^\alpha$, if there are sequences $(t^+_n)$ and $(t^-_n)$ with
\begin{equation}\label{rnowe}
t^\pm_n \to \pm\infty, \quad u(t^+_n)\to u_+ \quad \text{and} \quad u(t^-_n)\to u_- \aas n\to +\infty.
\end{equation}
Let $K\subset X^\alpha$ be a subset. We say that $K$ is \emph{invariant} with respect to $\Phi$, if for every $x\in K$ there is full solution $\sigma$ of $\Phi$ such that $\sigma(0) = x$ and $\sigma(\R) \subset K$. If $N\subset X^\alpha$ then we define \emph{maximal invariant set} as \label{def-izol}
\begin{align*}
\inv(N):=\inv(N, \Phi):=\{x\in N \ | \ & \text{ there is a solution } \sigma:\R\to X^\alpha \text{ of } \Phi \\
& \text{ such that } \sigma(0) = x \text{ and }\sigma(\R)\subset N\}.
\end{align*}
A closed invariant set $K$ is called \emph{isolated}, if there is a closed set $N\subset X^\alpha$ such that $$K=\inv(N)\subset\inter N.$$ In this case $N$ is called \emph{isolating neighborhood} for $K$. A subset $N\subset X^\alpha$ is \emph{admissible} with respect to $\Phi$, if for every sequences $(x_n)$ in $X^\alpha$, $(t_n)$ in $[0,+\infty)$ such that $t_n \to +\infty$ when $n\to \infty$, the inclusion $$\Phi([0,t_n]\times \{x_n\})\subset N \ffor n\ge 1,$$ implies that the set $\{\Phi(t_n,x_n) \ | \ n\ge 1\}$ is relatively compact in $X^\alpha$.

Let $\{\Phi^s\}_{s\in[0,1]}$ be a family of semiflows where $\Phi^s:[0,+\infty)\times X^\alpha \to X^\alpha$, $s\in[0,1]$. We say that $N\subset X^\alpha$ is \emph{admissible} with respect to $\{\Phi^s\}_{s\in[0,1]}$, if for every sequences $s_n\in[0,1]$, $(x_n)$ in $X^\alpha$ and $(t_n)$ in $[0,+\infty)$ such that $t_n \to +\infty$ when $n\to \infty$, the inclusion $$\Phi_{s_n}([0,t_n]\times \{x_n\})\subset N \ffor n\ge 1,$$ implies that the set $\{\Phi^{s_n}_n(t_n,x_n) \ | \ n\ge 1\}$ is relatively compact in $X^\alpha$. Furthermore the family $\{\Phi^s\}_{s\in[0,1]}$ is continuous provided $\Phi^{s_n}(t_n,x_n) \to \Phi^{s_0}(t_0,x_0)$ as long as $s_n\to s_0$, $t_n\to t_0$ and $x_n \to x_0$ as $n\to +\infty$.

From now on we write $\mathcal{S}(X^\alpha)=\mathcal{S}(X^\alpha, \Phi)$ for a class of invariant sets admitting an admissible isolated neighborhood. A special case of isolated neighborhood is {\em an isolating block}. To define it assume that $B\subset X^\alpha$ is a closed set and let $x\in\partial B$. We say that $x$ is \emph{a strict egress point} (resp. \emph{strict ingress point}, resp. \emph{bounce off point}), if for any solution $\sigma:[-\delta_1, \delta_2) \to X^\alpha$, where $\delta_1 \ge 0$ and $\delta_2 > 0$, of the semiflow $\Phi^s$ such that $\sigma(0) = x$ the following holds:
\begin{enumerate}
\item[(a)] there is $\varepsilon_2\in(0,\delta_2]$ such that $\sigma(t)\notin B$ (resp. $\sigma(t)\in \inter B$, resp. $\sigma(t)\notin B$) for $t\in(0,\varepsilon_2]$;
\item[(b)] if $\delta_1 > 0$ then there is $\varepsilon_1\in(0,\delta_1)$ such that $\sigma(t)\in \inter B$ (resp. $\sigma(t)\notin B$, resp. $\sigma(t)\notin B$) for $t\in[-\varepsilon_1, 0)$.
\end{enumerate}
We write $B^e$, $B^i$ and $B^b$ for the sets of strict egress points, strict ingress points and bounce off points, respectively. Furthermore, put $B^-:=B^e \cup B^b$. A close set $B\subset X^\alpha$ is \emph{an isolating block}, if $\partial B = B^e\cup B^i\cup B^e$ and the set $B^-$ is closed.

For $K\in \mathcal{S}(X^\alpha)$ we define the Conley index (\emph{the homotopy index}) of $K$ as the homotopy type $h(\Phi^s, K)$ of pointed space given by
\begin{equation*}\label{def-con}
h(\Phi^s, K):=\begin{cases}[B/B^-, [B^-]] & \text{ if } B^-\neq\emptyset; \\ [B\dot\cup \{c\}, c] & \text{ if } B^- = \emptyset\end{cases}
\end{equation*}
where $B/B^-$ is the quotient space and $B\dot\cup \{c\}$ is a disjoint sum of $B$ and the one point space $\{c\}$. It is known that the homotopy index is independent from the choice of isolating block of $K$ has the following properties: \\[5pt]
\makebox[8mm][l]{(H1)}\parbox[t]{118mm}{\emph{(Existence)} If $K\in\mathcal{S}(X^\alpha)$ and $h(\Phi, K)\neq \o 0$, then $K\neq\emptyset$.}\\[5pt]
\makebox[8mm][l]{(H2)}\parbox[t]{118mm}{\emph{(Addition)} If $K_1,K_2\in\mathcal{S}(X^\alpha)$ are disjoint sets, then $K_1\cup K_2\in\mathcal{S}(X^\alpha)$ and $$h(\Phi,K_1\cup K_2) = h(\Phi,K_1)\vee h(\Phi, K_2).$$}\\
\makebox[8mm][l]{(H3)}\parbox[t]{118mm}{\emph{(Multiplication)} Let $\Phi^i:[0,+\infty)\times X^\alpha \to X^\alpha$, for $i=1,2$ be semiflows. If $K_1,K_2\in\mathcal{S}(X^\alpha)$ then $K_1\times K_2 \in \mathcal{S}(\Phi^1\times\Phi^2, X^\alpha\times X^\alpha)$ and $$h(\Phi^1\times\Phi^2, K_1\times K_2) = h(\Phi^1, K_1)\wedge h(\Phi^2, K_2).$$}\\
\makebox[8mm][l]{(H4)}\parbox[t]{118mm}{\emph{(Homotopy invariance)} Let $\{\Phi^s:[0,+\infty)\times X^\alpha \to X^\alpha\}_{s\in[0,1]}$ be a continuous family of semiflows and let the set $N\subset X$ be admissible with respect to  this family. If for any $s\in[0,1]$ the set $N$ is an isolating neighborhood of $K_s:=\inv(\Phi^s, N)$, then $K_s\in\mathcal{S}(\Phi^s, X^\alpha)$ for $s\in[0,1]$ and $$h(\Phi^0, \inv(\Phi^0, N)) = h(\Phi^1, \inv(\Phi^1, N)).$$}

\section{Index formula for bounded orbits}

We will study the problem of existence of orbits connecting stationary points for the equations of the form
\begin{equation}\label{row-a-f}
\dot u(t) = - A u(t) + \lambda u(t) + F (u(t)), \qquad  t > 0,
\end{equation}
where $\lambda$ is an eigenvalue of a sectorial operator $A:X\supset D(A)\to X$ on a real Banach space $X$ with norm $\|\cdot\|$ and $F\colon X^\alpha \to X$ is a continuous map. Here $X^\alpha$ for $\alpha\in(0,1)$, is a fractional space given by $X^\alpha:=D((A + \delta I)^\alpha)$, where $\delta > 0$ is such that $A + \delta I$ is a positively defined operator. We assume that \\[5pt]
\noindent\makebox[9mm][l]{$(A1)$}\parbox[t][][t]{118mm}{$A$ has compact resolvent,}\\[5pt]
\noindent\makebox[9mm][l]{$(A2)$}\parbox[t][][t]{118mm}{there exists a Hilbert space $H$ equipped with a scalar product $\<\,\cdot\,, \,\cdot\,\>_H$ and norm $\|\cdot\|_H$ with a continuous injection $i:X \hookrightarrow H$,}\\[5pt] 
\noindent\makebox[9mm][l]{$(A3)$}\parbox[t][][t]{118mm}{there exists a linear self-adjoint operator $\h A:H\supset D(\h A) \to H$ with the property that $\mathrm{Gr}\,(A)\subset \mathrm{Gr}\,(\h A)$, where the inclusion is understood by the map $X \times X \xrightarrow{i\times i} H\times H$.}\\[5pt]   
\noindent\makebox[9mm][l]{$(F1)$}\parbox[t]{118mm}{given $x\in X^\alpha$ there is an open $V\subset X^\alpha$ of $x$ and $L > 0$ such that$$\|F(x_1) - F(x_2)\|\le L \|x_1 - x_2\|_\alpha\qquad\text{for}\quad x_1,x_2\in V;$$}\\
\noindent\makebox[9mm][l]{$(F2)$}\parbox[t][][t]{118mm}{there is a constant $m > 0$ such that $\|F(x)\|\le m$ for $x\in X^\alpha$.} \\[5pt] 
Let $\{S_A(t)\}_{t\ge 0}$ be a semigroup generated by $-A$ and let $J\subset\R$ be an interval. We say that a continuous map $u:J \to X^\alpha$ is \emph{a mild solution} of equation \eqref{row-a-f}, if
\begin{equation*}
u(t) = S_A(t - t')u(t') + \int_{t'}^t S_A(t - \tau)F(s,u(\tau)) \d \tau
\end{equation*}
for every $t,t' \in J$, $t' < t$. Since $A$ is sectorial and assumptions $(F1)$ and $(F2)$ hold, it is known (see \cite[Theorem 3.3.3, Corollary 3.3.5]{MR610244}) that for any $x\in X^\alpha$ equation \eqref{row-a-f} admits a unique mild solution $u(t; x):[0,+\infty)\to X^\alpha$ starting at $x$. Hence we define a semiflow $\Phi:[0,+\infty)\times X^\alpha \to X^\alpha$ associated with this equation by $$\Phi(t,x):= u(t;x) \ffor t\ge 0, \ \ x\in X^\alpha.$$
In view of Theorems 3.4 and 3.5 from \cite{Kok2}, the semiflow $\Phi$ is continuous and admissible with respect to any bounded set $N\subset X^\alpha$.
\begin{remark}\label{rem-pom}
The spectrum $\sigma(A)$ consists of the sequence (possibly finite) of real eigenvalues. Indeed, the operator $A$ has compact resolvents which implies that $$\sigma(A) = \{\lambda_i \ | \ i\ge 1\}\subset\C$$ and this set is finite or $|\lambda_i|\to +\infty$ when $n\to +\infty$. Furthermore, if $\lambda\in\C$ is a complex eigenvalue of $A$, then, by $(A3)$, it is also a complex eigenvalue of the symmetric operator $\h A$ and hence $\lambda$ is a real number. \hfill $\square$
\end{remark}
From the above remark it follows that the spectrum $\sigma(A)$ of the operator $A$ consists of the sequence of eigenvalues $$\lambda_1 < \lambda_2 < \ldots < \lambda_i < \lambda_{i+1} < \ldots$$ which is finite or $\lambda_i \to +\infty$ when $i\to +\infty$. Furthermore, by Theorem 2.1 from \cite{Kok2} we have the following
\begin{theorem}\label{th:10ab}
Under the assumptions $(A1)$, $(A2)$ and $(A3)$, if $\lambda = \lambda_k$ for some $k\ge 1$, is an eigenvalue of $A$, then there exists a a decomposition $X = X_+\oplus X_-\oplus X_0$ on closed subspaces, such that 
\begin{equation}\label{inkl-sem}
S_A(t)X_i \subset X_i \ffor t\ge 0, \ i\in\{0,-,+\},
\end{equation}
and the following holds: 
\begin{enumerate}
\item[(i)] $X_0 = \ker (\lambda I - A)$ and the space $X_-$ is such that 
\begin{equation*} 
   X_-= \{0\} \text{ provided } k=1 \  \text{ and } \   X_-=\bigoplus_{i=1}^{k-1} \ker(\lambda_i I - A) \text{ provided } k\ge 2.
\end{equation*}
Therefore $\dim X_- = 0$ for $k=1$ and $\dim X_- = \sum_{i=1}^{k-1} \dim\ker(\lambda_i I - A)$ for $k\ge 2$.
\item[(ii)] there are $c, M > 0$ with the property that 
\begin{align}\label{ine11}
     & \|A_\delta^\alpha S_A(t)x\| \le M e^{- (\lambda + c) t} t^{-\alpha} \|x\| &     \text{if} \quad t > 0, \ x\in X_+,\\ \label{ine22}
     & \|e^{\lambda t}S_A(t)x\| \le M e^{- c t}\|x\| & \text{if} \quad t\ge 0, \ x\in X_+, \\ \label{ine33}
     & \|e^{\lambda t}S_A(t)x\| \le M e^{c t}\|x\| &\text{if} \quad t\le 0, \ x\in X_-,
\end{align} 
\item[(iii)] $X_0$, $X_-$ and $X_+$ are mutually orthogonal spaces, that is, 
    \begin{equation}\label{ort}
    \<i (u_l),i(u_m)\>_H = 0 \ffor u_l\in X_l \ \  \text{and} \ \ u_m\in X_m
    \end{equation}
    where $l,m\in\{0,-,+\}$, $l\neq m$.
\end{enumerate}
\end{theorem}
Assume that $P, Q_\pm:X\to X$ are projections given 
\begin{equation}\label{wzz1ab}
    P x = x_0 \ \text{ and } \ Q_\pm x = x_\pm\qquad\text{for}\quad x\in X 
\end{equation}
where $x = x_+ + x_0 + x_-$ if $x_i\in X_i$, $i\in \{0,-,+\}$. Put $Q:= Q_- + Q_+$. In view of the continuity of the inclusion  $X^\alpha\subset X$, we have $X^\alpha= X_0\oplus X^\alpha_-\oplus X^\alpha_+$, where $$X^\alpha_-:= X^\alpha\cap X_-, \ X^\alpha_+:=X^\alpha\cap X_+$$ are closed spaces. Hence $P$ and $Q_\pm$ can be seen as continuous maps $P, Q_\pm:X^\alpha\to X^\alpha$. Observe also that \eqref{inkl-sem} implies  
\begin{align}\label{ds1}
S_A(t)Px = PS_A(t)x \ \text{ and } \ S_A(t)Q_\pm x = Q_\pm S_A(t) x \ffor t\ge 0, \ x\in X.
\end{align}
\begin{remark}\label{rem-res}
If the equation \eqref{row-a-f} is at resonance at infinity the problem of existence of orbits connecting stationary points may not have solution for general nonlinearity $F$. To see this it is enough to take $F(x) = y_0$ for $x\in X$, where $y_0\in\ker (\lambda I - A)\setminus\{0\}$. Indeed, if $u:\R\to X^\alpha$ is an orbit connecting stationary points for equation \eqref{row-a-f}, then $$u(t) = e^{\lambda (t - t')}S_A(t - t')u(t') + \int_{t'}^t e^{\lambda (t - \tau)}S_A(t - \tau) y_0 \d\tau \ffor t > t'.$$ Since $\ker (\lambda I - A) \subset \ker(I - e^{\lambda t}S_A(t))$ for $t\ge 0$ it follows that $$u(t) = e^{\lambda (t - t')}S_A(t - t')u(t') + (t - t')y_0 \ffor t > t'.$$
Therefore by \eqref{ds1} $$Pu(t) = e^{\lambda (t - t')}S_A(t - t')Pu(t') + (t - t')Py_0 = Pu(t') + (t - t')y_0 \ffor t > t',$$ and hence $Pu(h) = Pu(0) + h y_0$ for $h\ge 0$. Let sequences $(t^+_n)$, $(t^-_n)$ and stationary points $u_+,u_-\in X^\alpha$ be such that \eqref{rnowe} holds. Putting $h:=t^+_n$ we obtain a contradiction because $y_0\neq 0$.
\end{remark}
To overcome the difficulty discussed in the above remark we shall introduce geometric conditions for $F$ which will guarantee the existence of bounded orbits for equation \eqref{row-a-f}: \label{g1g2}
\begin{equation*}\leqno{(G1)}
\quad\left\{\begin{aligned}
&  \text{for every ball $B\subset X^\alpha_+\oplus X^\alpha_-$ there is $R > 0$ such that} \\
&  \<F(x + y), x\>_H > 0 \ \text{ for } (x, y)\in X_0\times B  \ \text{ such that } \|x\|_H\ge R,
\end{aligned}\right.
\end{equation*}
\begin{equation*}\leqno{(G2)}
\quad\left\{\begin{aligned}
& \text{for every ball $B\subset X^\alpha_+\oplus X^\alpha_-$ there is $R > 0$ such that} \\
& \<F(x + y), x\>_H < 0 \ \text{ for } (x, y)\in X_0\times B \ \text{ such that } \|x\|_H\ge R.
\end{aligned}\right.
\end{equation*}
Now we proceed to the main result of this section, namely \emph{the index formula for bounded orbits}. It is a tool to determining the Conley index for the maximal invariant set contained in appropriately large ball in terms of geometrical conditions $(G1)$ and $(G2)$. This theorem can be used directly to prove the existence of bounded orbits for the equation \eqref{row-a-f} or can be applied to prove the existence of orbits connecting stationary points, which is studying in the subsequent theorem. Write
$$d_0:= 0 \quad \text{ and } \quad d_l:=\sum_{i=1}^l \dim\ker(\lambda_i I - A) \quad\text{ for } \quad l\ge 1.$$
\begin{theorem}\label{th-ind-orbi}
If $\lambda = \lambda_k$ for some $k\ge 1$, then there is a closed isolated neighborhood $N\subset X^\alpha$, admissible with respect to the semiflow $\Phi$ such that $0\in \inter N$ and, for $K:=\inv(N,\Phi)$, the following statements hold:
\begin{enumerate}
\item[(a)] if condition $(G1)$ is satisfied, then $h(\Phi, K) = \Sigma^{d_k}$, \\[-10pt]
\item[(b)] if condition $(G2)$ is satisfied, then $h(\Phi, K) = \Sigma^{d_{k-1}}$.
\end{enumerate}
\end{theorem}
In the proof of the above theorem we will consider the family of equations
\begin{equation}\label{A-G}
\dot u(t)  = - A u(t) + \lambda u(t) + G(s, u(t)), \qquad  t > 0.
\end{equation}
where $G:[0,1]\times X^\alpha \to X$ is a map given by
\begin{equation}
G(s, x) := PF(s Q x + P x) + s QF(s Q x + P x) \ffor (s, x)\in[0,1] \times X^\alpha.
\end{equation}

Since $F$ is locally lipschitz, it is not difficult to check that $G(s,\,\cdot\,)$ satisfies condition $(F1)$ for every $s\in [0,1]$. Furthermore by $(F2)$, for any $s\in[0,1]$ and $x\in X^\alpha$ we have
\begin{equation}
\begin{aligned}\label{rem-imp}
\|G(s, x)\| & = \|PF(s Q x + P x) + s QF(s Q x + P x)\| \\
& \le \|P\| \|F(s Q x + P x)\| + \|Q\| \|F(s Q x + P x)\| \\
& \le m ( \|P\|  + \|Q\| ):=m_0,
\end{aligned}
\end{equation}
which shows that condition $(F2)$ is satisfied. Therefore, for any $s\in[0,1]$, one can define a semiflow $\Psi^s:[0,+\infty)\times X^\alpha \to X^\alpha$ given by the formula $$\Psi^s(t, x):=u(t;s, x) \ffor t\in[0,+\infty), \ x\in X^\alpha,$$ where $u(\,\cdot\,;s,x):[0,+\infty)\to X^\alpha$ is a weak solution of \eqref{A-G} starting at $x$. Theorems 3.4 and 3.5 from \cite{Kok2} say that the family $\{\Psi^s\}_{s\in[0,1]}$ is continuous and any bonded set is admissible. Hence we have a homotopy between $\Psi^1 = \Phi$ and the semiflow $\Psi^0$ associated with
\begin{equation*}
\dot u(t)  = - A u(t) + \lambda u(t) + PF(Pu(t)), \qquad  t > 0.
\end{equation*}
Note that every solution $u:[0,+\infty) \to X^\alpha$ of this equation satisfies the formula
\begin{equation}
u(t) = e^{\lambda t}S_A(t)u(0) + \int_0^t PF(Pu(\tau)) \d \tau \ffor t \ge 0
\end{equation}
Let $\psi_1:[0,+\infty)\times X^\alpha_-\oplus X^\alpha_+ \to X^\alpha_-\oplus X^\alpha_+$ be a semiflow given by $$\psi_1(t,x):= e^{\lambda t}S_A(t) x \ffor t\in [0,+\infty), \ x\in X^\alpha_-\oplus X^\alpha_+$$
and let semiflow $\psi_2:[0,+\infty)\times X_0\to X_0$ be associated with the equation
\begin{equation*}
\dot u(t)  = PF(u(t)),    \qquad  t > 0.
\end{equation*}
Then it is easy to see that  $$\Psi^0(t, x) = \psi_1(t, Qx) + \psi_2(t, Px) \ffor t\in[0,+\infty), \ x\in X^\alpha,$$ and therefore the semiflow $\Psi^0$ is equivalent with the product of $\psi_1$ and $\psi_2$, that is, for any $t\ge 0$, and $(x,y)\in (X^\alpha_-\oplus X^\alpha_+) \times X_0$ we have
\begin{equation}\label{rown1}
    \Psi^0(t, U(x,y)) = U(\psi_1(t, x), \psi_2(t, y))
\end{equation}
where $U:(X^\alpha_-\oplus X^\alpha_+) \times X_0 \to X^\alpha$ is defined by $U(x,y) = x + y$ for $(x,y)\in (X^\alpha_-\oplus X^\alpha_+) \times X_0$. \\

In the first step we prove the following lemma, which provides some {\em a priori} bounds for solutions of the equation \eqref{A-G}.
\begin{lemma}\label{prop:3}
There is $R > 0$ such that for any $s\in[0,1]$ and for any bounded full solution $u=u_s:\R \to X^\alpha$ for the semiflow $\Psi^s$, the following inequality holds
\begin{equation*}
\|Qu(t)\|_\alpha\le R \ffor t\in\R.
\end{equation*}
\end{lemma}
\noindent\textbf{Proof.} Let $s\in[0,1]$ be fixed and let $u=u_s:\R \to X^\alpha$ be a full solution for the equation \eqref{A-G}. Since $Q_-,Q_+:X^\alpha\to X^\alpha$ are bounded operators, the sets $\{Q_+u(t) \ | \ t\le 0\}$ and $\{Q_-u(t) \ | \ t\ge 0\}$ are bounded in $X^\alpha$. We prove that
\begin{equation}\label{fs9}
\|Q_+u(t)\|_\alpha \le m_0 M\|Q_+\|_{L(X)}\left(e^{-c}/c  + \frac{1}{1 - \alpha}\right) =:R_1 \ffor t\in\R.
\end{equation}
Indeed, since $u$ is full solution, we have the equality $\Psi^s(t-t',u(t')) = u(t)$ for $t,t'\in\R$, $t\ge t'$, which implies that
\begin{equation}\label{ful-sol}
u(t) = e^{\lambda(t - t')}S_A(t - t')u(t') + \int_{t'}^{t} e^{\lambda(t - \tau)}S_A(t - \tau)G(s,u(\tau)) \d \tau
\end{equation}
for $t\ge t'$. By \eqref{ds1} one find that
\begin{equation}\label{fs2}
Q_+u(t) = e^{\lambda(t - t')}S_A(t - t')Q_+u(t') + \int_{t'}^t e^{\lambda(t - \tau)}S_A(t - \tau)Q_+G(s,u(\tau)) \d \tau
\end{equation}
for $t\ge t'$. Since $X^\alpha$ embeds continuously in $X$, there is a constant $C > 0$ such that $\|x\| \le C \|x\|_\alpha$ for $x\in X^\alpha$. Furthermore, by the inequality \eqref{ine11}, there are constants $c, M>0$ such that
\begin{align*}
& \|e^{\lambda(t - t')}S_A(t - t')Q_+u(t')\|_\alpha = \|A_\delta^\alpha e^{\lambda(t - t')}S_A(t - t')Q_+u(t')\| \\
& \qquad \qquad\le M \frac{e^{- c (t - t')}}{(t - t')^\alpha} \, \|Q_+u(t')\| \le C M \frac{e^{- c (t - t')}}{(t - t')^\alpha} \, \|Q_+u(t')\|_\alpha
\end{align*}
for $t,t'\in\R$, $t>t'$. Hence, the boundedness of $\{Q_+u(t) \ | \ t\le 0\}$ implies that
\begin{equation}\label{fs3}
\|e^{\lambda(t - t')}S_A(t - t')Q_+u(t')\|_\alpha \to 0 \aas t'\to - \infty.
\end{equation}
In view of the inequality \eqref{ine11}, we deduce that
\begin{equation}\label{row-now2}
\begin{aligned}
& \|Q_+u(t)\|_\alpha \le \|e^{\lambda(t - t')}S_A(t - t')Q_+u(t')\|_\alpha \\
& \hspace{30mm} + \int_{t'}^t \|A_\delta^\alpha e^{\lambda(t - \tau)}S_A(t - \tau)Q_+G(s,u(\tau))\| \d \tau \\
& \le \|e^{\lambda(t - t')}S_A(t - t')Q_+u(t')\|_\alpha + M \int_{t'}^t \frac{e^{- c (t - \tau)}}{(t - \tau)^\alpha} \, \|Q_+G(s,u(\tau))\| \d \tau \\
& \le \|e^{\lambda(t - t')}S_A(t - t')Q_+u(t')\|_\alpha + m_0 M \|Q_+\|_{L(X)} \int_{t'}^t  \frac{e^{- c (t - \tau)}}{(t - \tau)^\alpha} \d \tau.
\end{aligned}
\end{equation}
Furthermore, taking $t, t'\in\R$ such that $t-t' > 1$, one find that
\begin{equation*}
\begin{aligned}
& \int_{t'}^t \frac{e^{- c (t - \tau)} }{(t - \tau)^\alpha} \d \tau
= \int_{t'}^{t-1} \frac{e^{- c (t - \tau)}}{(t - \tau)^\alpha} \d \tau +
\int_{t - 1}^t \frac{e^{- c (t - \tau)}}{(t - \tau)^\alpha} \d \tau \\
& \qquad \qquad \le \int_{t'}^{t-1} e^{- c (t - \tau)} \d \tau +
\int_{t - 1}^t \frac{1}{(t - \tau)^\alpha} \d \tau = (e^{-c} - e^{c (t' - t)})/c + 1/(1 - \alpha),
\end{aligned}
\end{equation*}
and hence, by \eqref{row-now2}, it follows that
\begin{equation*}
\begin{aligned}
\|Q_+u(t)\|_\alpha & \le \|e^{\lambda(t - t')}S_A(t - t')Q_+u(t')\|_\alpha  \\
& \hspace{10mm} +  m_0 M \|Q_+\|_{L(X)} \left((e^{-c} - e^{c(t' - t)})/ c + \frac{1}{1 - \alpha}\right).
\end{aligned}
\end{equation*}
Using \eqref{fs3} and passing to limit with $t'\to -\infty$ we infer that \eqref{fs9} is satisfied.
Since $X_-$ is finite dimensional there is constant $C' > 0$ such that
\begin{equation}\label{norm-in}
    \|x\|_\alpha \le C'\|x\| \ffor x\in X_-.
\end{equation}
We show that the boundedness of $\{Q_-u(t) \ | \ t\ge 0\}$ in $X^\alpha$ implies
\begin{equation}\label{fs7}
\|Q_-u(t)\|_\alpha \le m_0 C'M \|Q_-\|_{L(X)}/c :=R_2 \ffor t\in\R.
\end{equation}
In order to get this inequality we apply the operator $Q_-$ on equation \eqref{ful-sol} and, by \eqref{ds1}, we obtain
\begin{equation*}
Q_-u(t) = e^{\lambda(t - t')}S_A(t - t')Q_-u(t') + \int_{t'}^t e^{\lambda(t - \tau)}S_A(t - \tau)Q_-G(s,u(\tau)) \d \tau.
\end{equation*}
In consequence, for $t,t'\in\R$ and $t\ge t'$, we have
\begin{equation}
e^{\lambda (t' - t)}S_A(t' - t) Q_-u(t) = Q_-u(t') + \int_{t'}^t e^{\lambda (t' - \tau)}
S_A(t' - \tau)Q_-G(s,u(\tau)) \d \tau,
\end{equation}
where we used the fact that the family $\{S_A(t)\}_{t\ge 0}$ extends on the space $X_-$ to the $C_0$ group of bounded operators. The inequality \eqref{ine33} implies that
\begin{align*}
\|Q_-u(t')\| & \le \|e^{\lambda (t' - t)}S_A(t' - t) Q_-u(t)\| \\
& \qquad + \int_{t'}^t \|e^{\lambda (t' - \tau)} S_A(t' - \tau)Q_-G(s,u(\tau))\| \d \tau \\
& \le \|e^{\lambda (t' - t)}S_A(t' - t) Q_-u(t)\| + \int_{t'}^t M e^{c (t' - \tau)}\|Q_-G(s,u(\tau))\| \d \tau \\
& \le M e^{c(t' - t)} \|Q_-u(t)\| + m_0 M \|Q_-\|_{L(X)} \int_{t'}^t e^{c (t' - \tau)} \d \tau \\
& \le C M e^{c(t' - t)} \|Q_-u(t)\|_\alpha + m_0 M \|Q_-\|_{L(X)} \left(1 - e^{c (t' - t)}\right)/c.
\end{align*}
Therefore after passing to the limit with $t\to +\infty$ we infer that
\begin{equation}
\|Q_-u(t')\| \le  m_0 M \|Q_-\| /c \ffor t'\in\R
\end{equation}
which along with \eqref{norm-in} gives \eqref{fs7}. Hence, combining \eqref{fs9} and \eqref{fs7} gives
\begin{align*}
\|Qu(t)\|_\alpha & \le \|Q_-u(t)\|_\alpha + \|Q_+u(t)\|_\alpha \le \|Q_-\|_{L(X^\alpha)}\|u(t)\|_\alpha + \|Q_+\|_{L(X^\alpha)}\|u(t)\|_\alpha  \\ & \le R_1 \|Q_+\|_{L(X^\alpha)} + R_2 \|Q_-\|_{L(X^\alpha)}:= R,
\end{align*}
for $t\in\R$, which completes the proof. \hfill $\square$ \\

\noindent\textbf{Proof of Theorem \ref{th-ind-orbi}.} \textbf{Step 1.} Proposition \ref{prop:3} says that there is a constant $R_1 > 0$ such that for any $s\in[0,1]$ and for any full solution $u=u_s:\R \to X^\alpha$ for the semiflow $\Psi^s$ which is bounded in $X^\alpha$ we have
\begin{equation}\label{rrow}
\|Qu(t)\|_\alpha\le R_1 \ffor t\in\R.
\end{equation}
Suppose that condition $(G1)$ is satisfied. If we take $N_1 := \{x\in X^\alpha_+\oplus X^\alpha_+ \ | \ \|x\|_\alpha \le R_1 + 1\}$, then by \eqref{ort}, there is $R_2 > 0$ such that
\begin{equation}\label{rrow2}
\<PF(x + y), x\>_H > 0 \quad\text{ for } \ \  (y, x)\in N_1 \times X_0 \ \text{ with } \|x\|_H\ge R_2.
\end{equation}
Similarly, if we assume that condition $(G2)$ is satisfied, then there is $R_2 > 0$ such that
\begin{equation}\label{rrow3}
\<PF(x+y), x\>_H < 0 \quad\text{ for } \ \  (y, x)\in N_1 \times X_0 \ \text{ with } \|x\|_H\ge R_2.
\end{equation}
Define $N_2 := \{x\in X_0 \ | \ \|x\|_H \le R_2\}$. \\[5pt]
\textbf{Step 2.} We claim that $N:= N_1\oplus N_2$ is an isolating neighborhood for the family $\{\Psi^s\}_{s\in[0,1]}$. To proof this, let $u:=u_s:\R\to X^\alpha$ be a full solution for the semiflow $\Psi^s$, where $s\in[0,1]$, such that $u(\R)\subset N$ and $u(\R)\cap\partial N \neq \emptyset$. Without loss of generality one can assume that $u(0)\in\partial N$. Then we have either $\|Qu(0)\|_\alpha = R_1 + 1$ and $\|Pu(0)\|_H \le R_2$ or $\|Qu(0)\|_\alpha \le R_1 + 1$ and $\|Pu(0)\|_H = R_2$.
From \eqref{rrow}, it follows that the later holds. Since $u$ satisfies the integral formula
\begin{equation*}
u(t) = e^{\lambda(t - t')}S_A(t - t')u(t') + \int_{t'}^{t} e^{\lambda(t - \tau)}S_A(t - \tau)G(s,u(\tau)) \d \tau \ffor t\ge t',
\end{equation*}
from \eqref{ds1} we obtain
\begin{equation}\label{f1}
Pu(t) = e^{\lambda(t - t')}S_A(t - t')Pu(t') + \int_{t'}^{t} e^{\lambda(t - \tau)}S_A(t - \tau)PG(s,u(\tau)) \d \tau.
\end{equation}
On the other hand $\ker (\lambda I - A)\subset \ker (I - e^{\lambda t}S_A(t))$ for $t\ge 0$, and therefore \eqref{f1} takes the form
\begin{equation*}
Pu(t) = Pu(t') + \int_{t'}^t PF(sQu(\tau) + Pu(\tau)) \d \tau.
\end{equation*}
Endowing $X_0$ with norm $\|\cdot\|_H$ we see that the map $\R\ni t \to Pu(t) \in X_0$ is continuously differentiable on $\R$ and furthermore
\begin{equation*}
\frac{d}{dt} \|Pu(t)\|_H^2 = 2 \left\< \frac{d}{dt}P u(t), Pu(t)\right\>_H = 2\<PF(sQu(t) + Pu(t)), Pu(t)\>_H
\end{equation*}
for $t\in\R$. Since $\|Pu(0)\|_H^2 = R_2$, combining \eqref{rrow}, \eqref{rrow2} and \eqref{rrow3} we infer that
$$\frac{d}{dt}\|Pu(t)\|_H^2|_{t=0} = 2\<PF(sQu(0) + Pu(0)), Pu(0)\>_H\neq 0$$
if either $(G1)$ or $(G2)$ is satisfied. Therefore we find that the set $\{Pu(t) \ | \ t\in\R\}$ is not contained in $N_2$, which contradicts the inclusion $\{u(t) \ | \ t\in\R\}\subset N$. \\[5pt]
\textbf{Step 3.} Now we verify that $B:=N_2$ is an isolating block for the semiflow $\psi_2$ and
\begin{equation}\label{przypadkiblok}
\begin{cases}
(B,B^-)=(N_2, \partial N_2) & \text{ if condition $(G1)$ is satisfied;} \\
(B,B^-)=(N_2, \emptyset) & \text{ if condition $(G2)$ is satisfied.}
\end{cases}
\end{equation}
Assume that $u:[-\delta_2,\delta_1) \to X_0$, where $\delta_1 > 0$, $\delta_2 \ge 0$, is a solution for $\psi_2$ such that $u(0)\in \partial N_2$. Then
$$u(t) = u(0) + \int_0^t PF(u(\tau)) \d \tau \ffor t\in [-\delta_2,\delta_1),$$
which implies that the map $[-\delta_2,\delta_1)\ni t \to u(t) \in X_0$ is continuously differentiable when the space $X_0$ is endowed with the norm $\|\cdot\|_H$ and furthermore
\begin{equation*}
\begin{aligned}
\frac{d}{dt}\|u(t)\|_H^2 & = 2\<\dot u(t), u(t)\>_H = 2\<PF(u(t)), u(t)\>_H \ffor t\in [-\delta_2,\delta_1).
\end{aligned}
\end{equation*}
Since $\|u(0)\|_H = R_2$, combining \eqref{rrow2} and \eqref{rrow3} gives
\begin{equation}\label{przypadki}
\begin{cases}
\frac{d}{dt}\|u(t)\|_H^2|_{t=0} > 0 & \text{ if condition $(G1)$ is satisfied;} \\
\frac{d}{dt}\|u(t)\|_H^2|_{t=0} < 0 & \text{ if condition $(G2)$ is satisfied.}
\end{cases}
\end{equation}
The first inequality implies that in the case of $(G1)$ the pair $(B,B^-) := (N_2, \partial N_2)$ is an isolating block for the semiflow $\psi_2$. Similarly, in the case of the condition  $(G2)$, the second inequality \eqref{przypadki} shows that the pair $(N_2, \emptyset)$ is an isolating block for the semiflow $\psi_2$. \\[5pt]
\textbf{Step 4.} Applying Step 1 and homotopy invariance of Conley index we infer that
\begin{equation}\label{h1}
h(\Phi, K) = h(\Psi^1, K_1) = h(\Psi^0, K_0)
\end{equation}
where $K_s:=\inv (N, \Psi^s)$ for $s\in\{0,1\}$. Further, combining \eqref{ine22}, \eqref{ine33} and the fact that $$(\delta I  + A)^\alpha S_A(t)x = S_A(t) (\delta I  + A)^\alpha x \ffor x\in X^\alpha,$$ we deduce that there are constants $c,M > 0$ such that
\begin{alignat*}{2}
\|e^{\lambda t}S_A(t)x\|_\alpha &\le M e^{- c t}\|x\|_\alpha && \ffor t\ge 0, \ x\in X_+, \\
\|e^{\lambda t}S_A(t)x\|_\alpha &\le M e^{c t}\|x\|_\alpha && \ffor t\le 0, \ x\in X_-.
\end{alignat*}
Hence \cite[Theorem 11.1]{MR910097} shows that $N_1$ is an isolating neighborhood for $\psi_1$ and $K_0^1:=\inv(\psi_1, N_1)=\{0\}$ with
\begin{equation}\label{con-deg}
    h(\psi_1, K_0^1) = \Sigma^{\dim X_-} = \Sigma^{d_{k-1}}.
\end{equation}
In view of Step 3, the set $N_2$ is an isolating block for $\psi_2$. Take $K_0^2:=\inv(\psi_2, N_2)$. Combining \eqref{rown1}, \eqref{h1} and multiplication property of Conley Index we infer that
\begin{equation*}
h(\Phi, K) = h(\Psi^0, K_0) = h(\psi_1, K_0^1)\wedge h(\psi_2, K_0^2).
\end{equation*}
Therefore, by \eqref{con-deg}, we have
\begin{equation}\label{h2}
h(\Phi, K) = \Sigma^{\dim X_-}\wedge h(\psi_2, K_0^2) = \Sigma^{d_{k-1}}\wedge h(\psi_2, K_0^2).
\end{equation}
In the case of condition $(G1)$, the set $N_2$ is a ball in the Hilbert space $(X_0,\|\cdot\|_H)$ and hence the pair $(B,B^-)$ is homeomorphic with the pair of topological spaces where the first is a unit ball in the Euclidean space and the second space is its boundary (to see this it is enough to take orthogonal base in $X_0$). In consequence $$h(\psi_2, K_0^2)=\Sigma^{\dim X_0}.$$ Substituting this in \eqref{h2} we deduce that
\begin{equation*}
h(\Phi, K) = \Sigma^{d_{k-1}}\wedge \Sigma^{\dim X_0} = \Sigma^{d_k}
\end{equation*}
which proves point $(i)$. In the case of condition $(G2)$, from \eqref{przypadkiblok} it follows that the pair $(N_2, \emptyset)$ is an isolating block for the semiflow  $\psi_2$, which yields $h(\psi_2, K_0^2)=\Sigma^0$. Hence, by \eqref{h2}, we infer that
\begin{equation*}
h(\Phi, K) = \Sigma^{d_{k-1}}\wedge \Sigma^0 = \Sigma^{d_{k-1}}
\end{equation*}
which completes the proof of point $(ii)$. \hfill $\square$ \\

Now we apply Theorem \ref{th-ind-orbi} to study the orbits connecting stationary points. Before we do this, we make the following additional assumption on $F$: \\[5pt]
\noindent\makebox[9mm][l]{$(F4)$}\parbox[t][][t]{119mm}{$F(0) = 0$, and the map $F$ is differentiable at $0$ and there is $\mu\in\R$ such that $$DF(0)[x] = \mu x \ffor x\in X^\alpha.$$}\\

\noindent From $(F4)$ it follows that $0$ is a stationary point for the semiflow $\Phi$, that is $\Phi(t,0) = 0$ for $t\ge 0$. The following theorem is a criterion determining the existence of orbits connecting stationary points for equation  \eqref{row-a-f}.
\begin{theorem}\label{th-con-inde}
Let $\lambda=\lambda_k$ for some $k\ge 1$ and assume that one of the following conditions is satisfied \\[5pt]
\makebox[7mm][r]{(i)} \parbox[t]{115mm}{$(G1)$ holds and $\lambda_l < \lambda + \nu < \lambda_{l+1}$ where $\lambda_l \neq \lambda$;}\\[5pt]
\makebox[7mm][r]{(ii)} \parbox[t]{115mm}{$(G1)$ holds and $\lambda + \nu < \lambda_1$; }\\[5pt]
\makebox[7mm][r]{(iii)} \parbox[t]{115mm}{$(G2)$ holds, $\lambda_{l-1} < \lambda + \nu < \lambda_l$ and $\lambda \neq \lambda_l$, where $l\ge 2$;}\\[5pt]
\makebox[7mm][r]{(iv)} \parbox[t]{115mm}{$(G2)$ hold, $\lambda + \nu < \lambda_1$ and $\lambda \neq \lambda_1$.}\\[5pt]
Then there is a full nonzero solution $\sigma:\R \to X^\alpha$ for the semiflow $\Phi$ such that $$\lim_{t\to -\infty} u(t) = 0 \qquad\text{ or }\qquad \lim_{t\to +\infty} u(t) = 0.$$
\end{theorem}
In the proof of this theorem we need the following lemma.
\begin{lemma}{\em (see \cite[Theorem 3.5]{MR910097})}\label{con-ind-1}
If $\lambda + \mu \notin\sigma(A)$, then $\{0\}\in \mathcal{S}(\Phi, X^\alpha)$ and $h(\Phi,\{0\})= \Sigma^{b_l}$, where $b_l := 0$ if $\lambda + \mu < \lambda_1$ and $$b_l := \sum_{i=1}^l \dim\ker(\lambda_i I - A) \ \ \text{ if } \ \ \lambda_l < \lambda + \mu < \lambda_{l+1}.$$
\end{lemma}
\noindent\textbf{Proof of the Theorem \ref{th-con-inde}.} From Lemma \ref{con-ind-1}, it follows that $\{0\}\in \mathcal{S}(\Phi, X^\alpha)$ and $h(\Phi,K_0)= \Sigma^{b_l}$. Furthermore by Theorem \ref{th-ind-orbi} we infer that there is as isolating invariant set $K\in \mathcal{S}(\Phi, X^\alpha)$ such that $K_0\subset K$ and
\begin{equation}
h(\Phi, K) = \begin{cases}
\Sigma^{d_k} & \text{ if condition $(G1)$ is satisfied;}\\
\Sigma^{d_{k-1}} & \text{ if condition $(G2)$ is satisfied.}
\end{cases}
\end{equation}
Assumptions $(i)$--$(v)$ imply that $h(\Phi, K) \neq h(\Phi, K_0)$. Since $h(\Phi,K_0) \neq \o 0$, application of Theorems 11.6 and 11.5 from \cite{MR910097} completes the proof. \hfill $\square$ \\

\section{Applications}

Let $\Omega\subset\R^n$, $n\ge 1$, be an open bounded set such that the boundary $\partial\Omega$ is of class $C^\infty$. Let $\mathcal{A}$ be a differential operator with the Dirichlet boundary conditions given by: 
$$\mathcal{A} \x (x) = - \sum_{i,j=1}^n D_j(a_{ij}(x)D_i \x(x)) \ffor \x\in C^2(\o\Omega),$$ where $a_{ij} = a_{ji}\in C^2(\overline\Omega)$ for $1\le i,j\le n$ and 
$$\sum_{1\le i,j\le n}a_{ij}(x)\xi^i\xi^j \ge c_0 |\xi|^2 \ffor x\in\Omega, \ \xi\in\R^n,$$ where $c_0>0$. 
Assume that we have a continuous map $f:\Omega\times\R\times\R^n\to\R$ such that: \\[5pt] 
\noindent\makebox[22pt][l]{$(E1)$} \parbox[t][][t]{118mm}{there exists $L > 0$ such that for $x\in\Omega$, $s_1,s_2\in\R$ and $y_1,y_2\in\mathbb{R}^n$
    \begin{align*}
    |f(x,s_1,y_1) - f(x,s_2,y_2)| & \le L(|s_1 - s_2| + |y_1 - y_2|),
    \end{align*}}\\
\noindent\makebox[22pt][l]{$(E2)$} \parbox[t][][t]{118mm}{there exists a constant $m > 0$ such that $$|f(x,s,y)| \le m \ffor x\in\Omega, \ s\in\R, \ y\in\mathbb{R}^n.$$}\\[5pt]
Given $X:=L^p(\Omega)$, where $p\ge 1$, define the operator $A_p: X\supset D(A_p)\to X$ by
\begin{equation}
\begin{aligned}\label{op-a}
D(A_p) & := W^{2,p}_0(\Omega) := \mathrm{cl}_{W^{2,p}(\Omega)}\left\{\phi\in C^{2}(\overline \Omega) \ | \  \phi_{|\partial\Omega} = 0\right\}, \\
A_p \x & := \mathcal{A} \x \ffor \x\in D(A_p).
\end{aligned}
\end{equation}
\noindent It is well known (see e.g. \cite{MR1778284, MR500580}) that $A_p$ is positively defined sectorial operator on $X$. Let $X^\alpha := D(A_p^\alpha)$ for ($\alpha\in(0,1)$) be a fractional space with the norm
\begin{align*}
\|\x\|_\alpha := \|A_p^\alpha \x\| \ffor \x\in X^\alpha.
\end{align*}
From now on we assume that \\[5pt]
\noindent\makebox[10mm][l]{$(E3)$} \parbox[t][][t]{117mm}{$p\ge 2n$ and $\alpha\in(3/4,1)$.}
\begin{remark}\label{rem-pom2} 
$(a)$ Note that $A_p$ satisfies $(A1)$, $(A2)$ and $(A3)$. Since $A_p$ has compact resolvent (see e.g. \cite{MR1778284, MR500580}), the assumption $(A1)$ holds. Assume that $H:=L^2(\Omega)$ with the standard inner product and norm $$\<\x,\y\>_{L^2} := \int_\Omega \x(x) \y(x)\d x, \qquad\qquad \left(\int_\Omega |\x(x)|^2\d x\right)^{1/2} \ffor \x,\y\in H$$ and put $\h A:= A_2$. Then we see that the boundedness of $\Omega$ and the fact that $p\ge 2$ imply that there is a continuous embedding $i:L^p(\Omega) \hookrightarrow L^2(\Omega)$ and the assumption $(A2)$ is satisfied. Furthermore we have $D(A_p)\subset D(\h A)$ and $\h A \x = A_p \x$ and $\x \in D(A_p)$. This shows that $A_p \subset \h A$ in the sense of the inclusion $i\times i$. Since the operator $\h A$ is self-adjoint (see e.g. \cite{MR1778284}) we see that the assumption $(A3)$ is also satisfied. \\[5pt]
$(b)$ Remark \ref{rem-pom} shows that the spectrum $\sigma(A_p)$ of the operator $A_p$ consists of sequence of positive eigenvalues $$0 < \lambda_1 < \lambda_2 < \ldots < \lambda_i < \lambda_{i+1} < \ldots \ffor i\ge 1,$$ and furthermore $(\lambda_i)$ is finite or $\lambda_i \to +\infty$ when $i\to +\infty$. \\[5pt]
$(c)$ Note that the following inclusion is continuous
\begin{equation}\label{zan}
    X^\alpha\subset C^1(\o\Omega).
\end{equation}
Indeed, according to assumption $(E4)$ we have $\alpha\in (3/4, 1)$ and $p\ge 2n$, and hence $2\alpha - \frac{n}{p} > 1$. Therefore, the assertion i a consequence of \cite[Theorem 1.6.1]{MR610244}. \\[5pt]
$(d)$ If $1 \ge \alpha > \beta \ge 0$ then the inclusion $X^\alpha \subset X^\beta$ is continuous and compact as \cite[Theorem 1.4.8]{MR610244} says. \hfill $\square$
\end{remark}
According to the point $(c)$ of the above remark we can define a map $F\colon X^\alpha \to X$ given, for any $\x\in X^\alpha$, as
\begin{equation}\label{odwz-f}
    F(\x)(x) := f(x, \x(x), \nabla \x(x)) \ffor x\in\Omega.
\end{equation}
We call $F$ \emph{the Niemytzki operator} associated with $f$ and furthermore, it is easy to prove the following lemma
\begin{lemma}\label{lem-nem-op}
The map $F$ is well defined, continuous and satisfies assumption $(F1)$. Furthermore there is a constant $K > 0$ such that
\begin{equation}\label{sw}
\|F(\x)\| \le K \ffor \x\in X^\alpha.
\end{equation}
\end{lemma}

\subsection{Properties of Niemytzki operator}

We proceed to study conditions that should satisfy the mapping $f$ in order to the associated Niemytzki operator $F$ meets the introduced earlier geometrical conditions $(G1)$ and $(G2)$. We begin with the following theorem stating that the well-known \emph{Landesman-Lazer} conditions from \cite{MR0267269} are actually particular case of the geometrical conditions $(G1)$ and $(G2)$. 

\begin{theorem}\label{lem-est2}
Let $f_+,f_-\colon \Omega \to \mathbb{R}$ be continuous functions such that
\begin{align*}
f_+(x) = \lim_{s \to +\infty} f(x,s,y) \quad\text{and}\quad f_-(x) = \lim_{s \to -\infty} f(x,s,y)
\end{align*}
for $x\in\Omega$, uniformly for $y\in\R^n$ and let $B\subset X^\alpha_+\oplus X^\alpha_-$ be subset, bounded in the norm $\|\cdot\|_\alpha$. \\[5pt]
\makebox[7mm][l]{(i)}\parbox[t]{117mm}{Assume that \label{con-ll1-ll2}
$$\int_{\{u>0\}} f_+(x) \x(x) \,d x  + \int_{\{u<0\}} f_-(x) \x(x) \,d x > 0 \leqno{(LL1)}$$
for $\x\in\ker(\lambda I - A_p)\setminus\{0\}$. Then there exists $R > 0$ such that for any $(\z,\x)\in B \times X_0$ with $\|\x\|_{L^2} \ge R$, one has:
\begin{equation*}
\<F(\z + \x), \x\>_{L^2} > 0. 
\end{equation*}}\\
\makebox[7mm][l]{(ii)}\parbox[t]{117mm}{Assume that
$$\int_{\{u>0\}} f_+(x) \x(x) \,d x + \int_{\{u<0\}} f_-(x) \x(x) \,d x < 0 \leqno{(LL2)}$$
for $\x\in\ker(\lambda I - A_p)\setminus\{0\}$. Then there exists $R > 0$ such that for any $(\z,\x)\in B \times X_0$ with $\|\x\|_{L^2} \ge R$, one has: 
\begin{equation*}
\<F(\z + \x), \x\>_{L^2} < 0.
\end{equation*}}
\end{theorem}
\noindent\textbf{Proof.} The proof of is similar to the proof of Theorem 6.7 from \cite{Kok2}. \hfill $\square$ \\

In the subsequent lemma we show that the geometrical conditions $(G1)$ and $(G2)$ are also consequences of \emph{the strong resonance conditions} from \cite{MR713209}, \cite{MR597281}.

\begin{theorem}\label{lem-est3}
Assume that there is a continuous function $f_\infty \colon \Omega \to \mathbb{R}$, where $\Omega\subset\R^n$ ($n\ge 3$), such that
\begin{equation*}
f_\infty(x)  = \lim_{|s| \to +\infty} f(x,s,y)\cdot s
\end{equation*}
for $x\in\Omega$, uniformly for $y\in\R^n$ ad let $B\subset X^\alpha_+\oplus X^\alpha_-$ be a subset, bounded in the norm $\|\cdot\|_\alpha$. \\[5pt]
\makebox[7mm][l]{(i)}\parbox[t]{117mm}{If the following condition holds \label{con-sr1-sr2}
\begin{equation*}\leqno{(SR1)}
\quad\left\{\begin{aligned}\label{cond-sr3}
& \text{there is } h\in L^1(\Omega) \text{ such that } \\
& f(x,s,y)\cdot s \ge h(x) \text{ for } (x,s,y)\in\Omega\times\R\times\R^n
\text{ and } \\ & \int_\Omega f_\infty(x)\d x > 0,
\end{aligned}\right.
\end{equation*}
then there exists $R > 0$ such that for $(\y,\x)\in B \times X_0$ with $\|\x\|_{L^2}\ge R$
\begin{equation*}
\<F(\z + \x), \x\>_{L^2} > 0.
\end{equation*}}\\[5pt]
\makebox[7mm][l]{(ii)}\parbox[t]{117mm}{If the following condition holds
\begin{equation*}\leqno{(SR2)}
\quad\left\{\begin{aligned}\label{cond-sr4}
& \text{there is a function } h\in L^1(\Omega) \text{ such that } \\
& f(x,s,y)\cdot s \le h(x) \text{ for } (x,s,y)\in \Omega\times\R\times\R^n\text{ and } \\ & \int_\Omega f_\infty(x)\d x < 0,
\end{aligned}\right.
\end{equation*}
then there exists $R > 0$ such that for $(\z, \x)\in B\times X_0$ with $\|\x\|_{L^2}\ge R$
\begin{equation*}
\<F(\z + \x), \x\>_{L^2}< 0.
\end{equation*}}
\end{theorem}
 
\noindent\textbf{Proof.} The proof of is similar to the proof of Theorem 6.9 from \cite{Kok2}. \hfill $\square$ \\

\subsection{Existence of connecting orbits}

We shall consider parabolic equations of the form
\begin{equation}\label{A-eps-res-a}
u_t(t, x) = - \mathcal{A} \, u(t,x) + \lambda  u(t,x) + f(x, u(t,x), \nabla u(t,x)), \ \ t > 0, \ x\in\Omega.
\end{equation}
This equation may be written in the abstract form as
\begin{equation}\label{row-dr}
\dot u(t)  = - A_p u(t) + \lambda u(t) + F (u(t)), \qquad  t > 0.
\end{equation}
If $J\subset \R$ is an interval, then we say that $u:J\to X^\alpha$ is a solution of \eqref{A-eps-res-a}, if $u$ is a mild solution of \eqref{row-dr}. From Lemma \ref{lem-nem-op} it follows that $F$ satisfies $(F1)$, $(F2)$. Let $\Phi: [0,+\infty) \times X^\alpha \to X^\alpha$ be a semiflow associated with \eqref{row-dr} defined by $$\Phi(t,\x) := u(t; \x) \ffor t\in [0,+\infty), \ \x\in X^\alpha.$$
From now on we will also assume that \\[5pt]
\noindent\makebox[22pt][l]{$(E4)$} \parbox[t][][t]{115mm}{$f:\o\Omega\times \R\times\R^n\to\R$ is a map of class $C^1$ such that there is a constant $\nu\in\R$ such that $\nu=D_s f(x,0,0)$ for $x\in\Omega$. Furthermore $f(x,0,0) = 0$ and $D_y f(x,0,0) = 0$ for $x\in\Omega$.}
\begin{remark}\label{rem-con-f7}
By assumption $(E4)$ and \eqref{zan} and one can easily prove that that $F$ is differentiable at $0$ and its derivative $DF(0)\in L(X^\alpha, X)$ is of the form $$DF(0)[\x] = \nu \x \ffor \x\in X^\alpha,$$ and hence assumption $(F4)$ is satisfied.
\end{remark}
We start with the following {\em criterion with Landesman-Lazer conditions}
\begin{theorem}\label{th-crit-ogr}
Suppose that $f_+,f_-\colon \Omega \to \mathbb{R}$ are continuous functions such that
\begin{equation*}
f_+(x) = \lim_{s \to +\infty} f(x,s,y) \quad\text{and}\quad f_-(x) = \lim_{s \to -\infty} f(x,s,y)
\end{equation*}
for $x\in\Omega$, uniformly for $y\in\R^n$. If $\lambda=\lambda_k$ for some $k\ge 1$, then there is a full compact nonzero solution $u:\R \to X^\alpha$ for \eqref{A-eps-res-a} such that either $$\lim_{t\to -\infty} u(t) = 0 \qquad\text{ or }\qquad \lim_{t\to +\infty} u(t) = 0,$$ provided one of the following conditions is satisfied: \\[5pt]
\makebox[7mm][r]{(i)} \parbox[t]{115mm}{$(LL1)$ is satisfied holds and $\lambda_l < \lambda + \nu < \lambda_{l+1}$ where $\lambda_l \neq \lambda$;}\\[5pt]
\makebox[7mm][r]{(ii)} \parbox[t]{115mm}{$(LL1)$ is satisfied and $\lambda + \nu < \lambda_1$;} \\[5pt]
\makebox[7mm][r]{(iii)} \parbox[t]{115mm}{$(LL2)$ is satisfied, $\lambda_{l-1} < \lambda + \nu < \lambda_l$ and $\lambda \neq \lambda_l$, where $l\ge 2$;}\\[5pt]
\makebox[7mm][r]{(iv)} \parbox[t]{115mm}{$(LL2)$ is satisfied, $\lambda + \nu < \lambda_1$ and $\lambda \neq \lambda_1$.}
\end{theorem}
\noindent\textbf{Proof.} By Remark \ref{rem-pom2} $(a)$ and Remark \ref{rem-con-f7}, we deduce that assumptions $(A1)$, $(A2)$, $(A3)$ and $(F4)$ are satisfied. Furthermore, from Theorem \ref{lem-est2} it follows that, condition $(LL1)$ implies condition $(G1)$ and that condition $(LL2)$ implies condition $(G2)$. Therefore Theorem \ref{th-con-inde} completes the proof.  \hfill $\square$ \\

Now we proceed to the following \emph{criterion with strong resonance conditions}.

\begin{theorem}\label{th-crit-ogrsr}
Assume that $\Omega\subset\R^n$ where $n\ge 3$, is an open bounded set and let $f_\infty \colon \o\Omega \to \mathbb{R}$ be a continuous function such that
\begin{gather*}
f_\infty(x)  = \lim_{|s| \to +\infty} f(x,s,y)\cdot s \quad\text{for} \ x\in\Omega, \ \text{uniformly for }y\in\R^n.
\end{gather*} 
Given $\lambda=\lambda_k$ for $k\ge 1$, there is a full compact nonzero  solution $u:\R \to X^\alpha$ of equation \eqref{A-eps-res-a} such that either $$\lim_{t\to -\infty} u(t) = 0 \qquad\text{ or }\qquad \lim_{t\to +\infty} u(t) = 0,$$ provided one of the following conditions is satisfied:\\[5pt]
\makebox[7mm][r]{(i)} \parbox[t]{115mm}{$(SR1)$ is satisfied and $\lambda_l < \lambda + \nu < \lambda_{l+1}$ where $\lambda_l \neq \lambda$;}\\[5pt]
\makebox[7mm][r]{(ii)} \parbox[t]{115mm}{$(SR1)$ is satisfied and $\lambda + \nu < \lambda_1$;} \\[5pt]
\makebox[7mm][r]{(iii)} \parbox[t]{115mm}{$(SR2)$ is satisfied and $\lambda_{l-1} < \lambda + \nu < \lambda_l$ where $\lambda \neq \lambda_l$, $l\ge 2$;}\\[5pt]
\makebox[7mm][r]{(iv)} \parbox[t]{115mm}{$(SR2)$ is satisfied and $\lambda + \nu < \lambda_1$.}
\end{theorem}
\noindent\textbf{Proof.} Similarly as before, Remark \ref{rem-pom2} $(a)$ and Remark \ref{rem-con-f7} imply that assumptions $(A1)$, $(A2)$, $(A3)$ and $(F4)$ are satisfied. Furthermore, from Theorem \ref{lem-est3} it follows that condition $(SR1)$ implies $(G1)$ and that $(SR2)$ implies condition $(G2)$. Therefore the proof is completed after application of Theorem \ref{th-con-inde}.  \hfill $\square$ \\

\def\cprime{$'$} \def\polhk#1{\setbox0=\hbox{#1}{\ooalign{\hidewidth
  \lower1.5ex\hbox{`}\hidewidth\crcr\unhbox0}}} \def\cprime{$'$}
  \def\cprime{$'$} \def\cprime{$'$}
\providecommand{\bysame}{\leavevmode\hbox to3em{\hrulefill}\thinspace}
\providecommand{\MR}{\relax\ifhmode\unskip\space\fi MR }
% \MRhref is called by the amsart/book/proc definition of \MR.
\providecommand{\MRhref}[2]{%
  \href{http://www.ams.org/mathscinet-getitem?mr=#1}{#2}
}
\providecommand{\href}[2]{#2}

\parindent = 0 pt

\end{document}